\documentstyle[12pt]{article}

\input amssym.def
\input amssym.tex

\newcommand{\ga}{\alpha}     
\newcommand{\gb}{\beta}      
\renewcommand{\gg}{\gamma}   
\newcommand{\gd}{\delta}     
  
\newcommand{\gz}{\zeta}      
\newcommand{\gee}{\eta}      
\newcommand{\gth}{\theta}    
      
\newcommand{\gk}{\kappa}  
\newcommand{\gl}{\lambda}    
\newcommand{\gm}{\mu}        
\newcommand{\gn}{\nu}

\newcommand{\gr}{\rho}       
\newcommand{\gs}{\sigma}     
\newcommand{\gt}{\tau}

\newcommand{\ha}{\aleph}

\newcommand{\setof}[2]{{\{\; #1 \; \vert \; #2 \; \} } }
\newcommand{\seq}[1]{{\langle #1 \rangle} }
\newcommand{\card}[1]{{\vert #1 \vert} }
\newcommand{\ot}[1]{\hbox{o.t.($#1$)}}

\renewcommand{\models}{\vDash}

\newcommand{\dom}{{\rm dom}}
\newcommand{\rge}{{\rm rge}}
\newcommand{\crit}{{\rm crit}}

\newcommand{\cf}{{\rm cf}}
\newcommand{\lh}{{\rm lh}}

\newcommand{\lra}{\longrightarrow}

\newtheorem{definition}{Definition}%[chapter]
\newtheorem{theorem}{Theorem}%[chapter]
\newenvironment{proof}{\noindent{\bf Proof:}}{\nopagebreak\mbox{}\newline
 \makebox[\textwidth]{\hfill$\blacklozenge$}\par\bigskip}

\newcommand{\implies}{\Longrightarrow}

\catcode`\@=11
\def\@begintheorem#1#2{\rm \trivlist \item[\hskip \labelsep{\bf #1\ #2:}]}
\def\@opargbegintheorem#1#2#3{\rm \trivlist
      \item[\hskip \labelsep{\bf #1\ #2\ (#3):}]}
\catcode`\@=12

\newtheorem{lemma}{Lemma}%[chapter]

\newcommand{\gkp}{{\gk^+}}

\newcommand{\FP}{{\Bbb P}}
\newcommand{\FQ}{{\Bbb Q}}
\newcommand{\FR}{{\Bbb R}}
\newcommand{\FS}{{\Bbb S}}

\title{A consistency result on weak reflection}

\author{James Cummings
        \thanks{Supported by a Postdoctoral Fellowship at the Hebrew University.} \\
        {\tt \small cummings@math.huji.ac.il} \\
        Mirna D\v zamonja
        \thanks{Partially supported by the Basic Research Fund of
        the Israel Academy of Science, and a Postdoctoral Fellowship from
        the Hebrew University.} \\
        {\tt \small dzamonja@math.huji.ac.il} \\
        Saharon Shelah
        \thanks{Partially supported by the Basic Research Fund of
        the Israel Academy of Science. Paper number 571.} \\
        {\tt \small shelah@math.huji.ac.il} \\
        Hebrew University of Jerusalem}

\begin{document}

\baselineskip=16pt
\binoppenalty=10000
\relpenalty=10000

\maketitle

\begin{abstract}
   In this paper we study the notion of {\em strong non-reflection,}
 and its contrapositive {\em weak reflection.} We say
  {\em $\gth$ strongly non-reflects at $\gl$\/} iff there
 is a function
 $F: \gth \lra \gl$ such that for all $\ga < \gth$ with $\cf(\ga) = \gl$ there
 is $C$ club in $\ga$ such that $F \restriction C$ is strictly increasing.
 We prove that it is consistent to have a cardinal $\gth$ such that
 strong non-reflection and weak reflection each hold on an unbounded
 set of cardinals less than $\gth$.
 \footnote{The research for this paper was done in the period
 July 1994 -- January 1995.}

\end{abstract}

\newcommand{\Ref}{{\rm Ref}}
\newcommand{\SNR}{{\rm SNR}}

\section{Introduction}

   In this paper we study the notion of {\em strong non-reflection\/}, which was introduced
  in \cite{545} and is further studied in \cite{vaporware}.
 We prove that for a fixed $\gth$ we can have an unbounded set of
 cofinalities at which strong non-reflection holds, and an unbounded set where it
 fails.

\begin{definition} Let  $\gth$ be a  regular cardinal, and let $\gl$ be an ordinal with
 $\gl \ge \gth$.
\begin{itemize}
\item $S^\gl_\gth = \setof{\ga < \gl}{\cf(\ga) = \gth}$.
\item $S^\gl_{<\gth} = \setof{\ga < \gl}{\cf(\ga) < \gth}$. 
\item If $\gth$ is uncountable, then 
  {\em $\gl$ strongly non-reflects at $\gth$\/} iff there
 is a function
 $F: \gl \lra \gth$ such that for all $\ga \in S^\gl_\gth$ there
 is $C$ club in $\ga$ such that $F \restriction C$ is strictly increasing.
 We will write $\SNR(\gl,\gth)$ for this.
\item $\gl$ {\em weakly reflects at $\gth$\/} iff $\gl$ does not strongly
 non-reflect at $\gth$.
\end{itemize}
\end{definition}

  In \cite{545} D\v zamonja and Shelah prove some theorems connecting weak club
 principles, saturated ideals, and the ideal ${\cal I}[\gl,\gth)$ consisting
 of those $A \subseteq \gl$ such that there is $h: \gl \lra \gth$ increasing
 on a club at every point of $A \cap S^\gl_\gth$. In particular Theorem 2.5
 of that paper shows that a certain weak club principle is incompatible with
 saturated ideals at successors of singulars, and Theorem 2.8 connects weak
 reflection and the weak club principle.
%
%  Mirna --- Should we say more here? James
%

 In \cite{vaporware} strong non-reflection is used as a tool to show that different
 instances of stationary reflection are independent. For example it is shown there
 that ``every stationary subset of $S^{\ha_3}_{\ha_0}$ reflects at a point in $S^{\ha_3}_{\ha_2}$'' 
 is consistent with ``every stationary subset of $S^{\ha_3}_{\ha_1}$ has a non-reflecting
 stationary subset''.

   We  make a few  remarks about the definitions.
 The next lemma
 is implicit in Observation 1.2.3 from \cite{545}.

\begin{lemma} \label{lemma1} Let $\gl$ be an ordinal, and $\seq{\gl_i: i < \gl}$ a
 $\gl$-sequence of ordinals with $\cf(\gl_i) \neq \gth$. Let $F: \gl \lra \gth$
 and $F_i: \gl_i \lra \gth$
 witness strong non-reflection at $\gth$ for $\gl$ and each of the $\gl_i$.
 Then there is $G: \sum\limits_{i < \gl} \gl_i \lra \gth$ witnessing strong
 non-reflection for  $\sum\limits_{i < \gl} \gl_i$.
\end{lemma}

\begin{proof} Define $G(\sum\limits_{i<j} \gl_i) = F(j)$ and
 $G(\sum\limits_{i<j} \gl_i + \gn) = F_j(\gn)$ for $0 < \gn < \gl_j$.
 It is easy to check that this works.
\end{proof}

  It is proved in \cite{545}, using the previous lemma,  that the least $\gl$ which weakly reflects
 at $\gth$ is a regular cardinal greater than $\gth$.

  As the terminology suggests, there is a connection between weak reflection
 and the more familiar notion of {\em stationary reflection.}

\begin{definition} Let $\gk < \gm < \gn$ be  regular cardinals. Then 
 $\Ref(\gn, \gm, \gk)$ iff for every stationary $S \subseteq S^\gn_\gk$
 there is $\ga \in S^\gn_\gm$ such that $S \cap \ga$ is stationary in $\ga$.
 We will use also  $\Ref(\gn, \gm, <\gk)$ as a shorthand for
 ``$\forall \gd < \gk \; \Ref(\gn, \gm, \gd)$''.
\end{definition}

  The next fact shows that strong non-reflection at $\gth$ is antithetical to stationary
 reflection to points of cofinality $\gth$.

\begin{lemma} \label{motivation} Suppose that $\gl$ strongly non-reflects at $\gth$. Then for
 every stationary $S \subseteq \gl$ there is $T \subseteq S$ stationary
 such that $T \cap \ga$ is non-stationary for all $\ga \in S^\gl_\gth$.
\end{lemma}

\begin{proof} Use Fodor's Lemma to find $T$ on which $F: \gl \lra \gth$
 witnessing the strong non-reflection is constant. If $C \subseteq \ga$
 is a club on which $F$ is strictly increasing then $C$ meets $T$
 in at most one point.
\end{proof}

  It is not hard to
 see that if $\seq{C_\ga: \ga < \gk^+}$ is a $\Box_\gk$-sequence
 then the function $F: \ga \longmapsto \ot{C_\ga}$  witnesses that
 $\gkp$ strongly non-reflects at $\gk$. More is true, see 1.7 of \cite{545}.
 The following remark is immediate from the definition.

\begin{lemma} \label{upwards} If $\gth = \cf(\gth) < \gl < \gl^*$, and $\gl$ weakly reflects
 at $\gth$, then $\gl^*$ weakly reflects at $\gth$.
\end{lemma}

  We are now ready to state the main result. 

\begin{theorem}
\label{mainthm}
 Let GCH hold, let $\gth$ be regular,  and suppose that
 there are cardinals 
 $\seq{\gl_i, \gth_i, \gk_i : i < \gth}$
  such that for all $i < \gth$
\begin{enumerate}
\item $\gl_i = \cf(\gl_i) < \gl_i^+ < \gth_i = \cf(\gth_i) < \gk_i$, and $\gk_i$
 is measurable.
\item $\gl_i > (\sup\setof{\gk_j}{j<i})^{++}$.
 
\end{enumerate}

 {\em Then\/} there is a generic extension in which  $SNR(\gth,\gl_i)$
 and $\Ref(\gk_i, \gth_i, <\gth_i)$ for every $i$. In particular,
 by Lemmas \ref{motivation} and \ref{upwards}, $\gth$ weakly reflects
 at $\gth_i$. 
\end{theorem}

  The proof will involve two stages.
  First we force functions that witness the strong non-reflection
 at the points $\gl_i$, via
 an iterated forcing with Easton support. An important feature of the construction is
 that the individual steps in the forcing have an increasing degree of strategic closure,
 and at any stage a witness to the desired degree of strategic closure was added by the
 previous stages.

 We will show that the first stage preserves the measurability of all the $\gk_i$.
 In the second stage we will force with a product of the L\'evy collapses $Coll(\gth_i, < \gk_i)$,
 and use Baumgartner's argument from \cite{JEB2} to show that $\Ref(\gk_i, \gth_i, <\gth_i)$
 holds in the extension.

\section{Forcing strong non-reflection}

 Let $\gs$ and $\gl$ be regular cardinals with $\gs < \gl$.  
In this section we define a forcing $\FP(\gs,\gl)$ which adds a function from
 $\gl$ to $\gs$ witnessing strong non-reflection for $\gl$ at $\gs$. 
We could make the same definition for $\gl$ an arbitrary ordinal
greater than $\gs$, but for our purposes it will suffice to restrict
 ourselves to regular cardinals.

\begin{definition} Conditions in $\FP(\gs, \gl)$ are functions $p$ such that
 $\dom(p) < \gl$, $\rge(p) \subseteq \gs$, and for every $\gb \le \dom(p)$
 if $\gb \in S^\gl_\gs$ then there is is a club $C$ in $\gb$ such that
 $p \restriction C$ is strictly increasing.

   The condition $p$ extends the condition $q$ iff $p \supseteq q$. We write
   this as $p \le q$.
\end{definition}

   Clearly this forcing has at most $\gs^{<\gl}$ conditions, so enjoys the
 $(\gs^{<\gl})^+$-c.c. There are several pieces of information about
 the closure properties of the forcing that we will need later.

\begin{lemma} Let $\seq{p_\ga:\ga < \gr}$ be a strictly
 decreasing sequence of conditions in $\FP(\gs, \gl)$, where $\gr \in S^\gl_{\neq \gs}$.
 Then $p = \bigcup_\ga p_\ga$ is the greatest lower bound for the sequence.
\end{lemma}

\begin{proof} Notice that $\cf(\dom(p)) = \cf(\gr) \neq \gs$, so that
 if $\gb \le \dom(p)$ and $\cf(\gb) = \gs$ then $\gb \in \dom(p_\ga)$
 for some $\ga < \gr$. There is $C$ club in $\gb$ such that
 $p \restriction C = p_\ga \restriction C$ is strictly increasing.
\end{proof}

 We remind the reader of the notion of {\em strategic closure.}

\begin{definition} Let $\FP$ be a partial ordering, and let $\gee$
 be an ordinal.

\begin{enumerate}

\item
 The game $G(\FP, \gee)$ is played by two players
 I and II, who take turns to play elements $p_\ga$ of $\FP$
 for $ 0< \ga < \gee$, with player I playing at odd stages
 and player II at even stages (Nota bene: limit ordinals are even).

  The rules of the game are that the sequence that is played
 must be decreasing (not necessarily strictly decreasing),
 the first player who cannot make a move loses, and player II
 wins if play proceeds for $\gee$ stages.

\item
 $\FP$ is {\em $\gee$-strategically closed\/} iff 
 player II has a winning strategy in $G(\FP, \gee)$.

\item
 $\FP$ is {\em $<\gee$-strategically closed\/} iff for all $\gz < \gee$
 $\FP$ is  $\gz$-strategically closed.
 
\end{enumerate}

\end{definition}

 We say that a forcing notion $\FP$ is {\em $<\gl$-distributive\/} iff
 it does not add any $<\gl$-sequence of ordinals to the ground model
 (equivalently, the intersection of fewer than $\gl$ dense open sets
  is nonempty). The following lemma is easy.

\begin{lemma} If $\FP$ is $<\gl$-closed it is $\gl$-strategically closed,
 and if $\FP$ is $<\gl$-strategically closed it is $<\gl$-distributive.
\end{lemma}

 Notice that $\FP(\gs, \gl)$ will only contain conditions of lengths
 unbounded in $\gl$ if $\SNR(\gm, \gs)$ holds for all $\gm < \gl$.
 This condition is actually enough to make $\FP(\gs,\gl)$ be
 $<\gl$-strategically closed.

\begin{lemma} \label{stratclos1}
 Suppose that all $\gm \in [\gs, \gl)$ are strongly non-reflecting
 at $\gs$. {\em Then\/} $\FP(\gs, \gl)$ is  $<\gl$-strategically closed.
\end{lemma}

\begin{proof} Let $\gee < \gl$. If $\gee < \gs$ then player II can win with
 the following strategy;  he plays
 $p_{2\gg} =_{def} \bigcup_{\ga < 2\gg} p_\ga$.

  If $\gs \le \gee < \gl$ then by hypothesis there is a function $F: \gee \lra \gs$
 witnessing strong non-reflection. Player II will play 
 $p_{2\gg} =_{\rm def} (\bigcup_{\ga < 2\gg} p_\ga) \frown F(\gg)$.
 We check that this is a winning strategy.

 Let $2 \gd$ be an even stage of cofinality $\gs$ in $G(\FP(\gs, \gl), \gee)$.
There is $D$ club in $\gd$ such that $F \restriction D$ is strictly increasing.
 If we define $C = \setof{ \lh(p_{2\gg})}{\gg \in D}$ then $C$ witnesses
 that II does not lose at stage $2 \gd$.
\end{proof}

  We will be interested in forcing strong non-reflection to several values
  of $\gs$ simultaneously. For this we will use a certain dense subset
  of the $<\gl$-support product of the appropriate $\FP(\gs, \gl)$.

\begin{definition} Let $A \subseteq {\rm REG} \cap \gl$.
 Then $\FP(A, \gl)$ is the set of functions $p$ such that
\begin{enumerate}
\item $\dom(p) = (A \cap \gg) \times \gg$ for some $\gg < \gl$.
\item If $\dom(p) = (A \cap \gg) \times \gg$ and $\gs \in A \cap \gg$
 then $\ga < \gg \longmapsto p(\gs, \ga)$ is a condition in $\FP(\gs, \gl)$.
\end{enumerate}
  If $p, q \in \FP(A, \gl)$ then $p \le q$ iff $p$ extends $q$.
\end{definition}

 Clearly $\card{\FP(A, \gl)} \le \gl^{<\gl}$, so the forcing has the
 $(\gl^{<\gl})^+$-c.c. 
 We also record some information about the closure of the forcing.

\begin{lemma} Let $\seq{p_\ga:\ga < \gr}$ be a strictly
 decreasing sequence of conditions in $\FP(A, \gl)$, where $\cf(\gr) \notin A$.
Then the condition $p$ given by $p(i) =_{\rm def} \bigcup_{\ga < \gr} p_\ga(i)$
 is the greatest lower bound for the sequence. 
\end{lemma}

   The next lemma is easy, with a proof almost identical to that of
 Lemma \ref{stratclos1}.

\begin{lemma} \label{stratclos}
 \label{strat} Let $A$ and $\gl$ be as above. Suppose that for all $\gs \in A$,
  all $\gm \in [\gs, \gl)$ are strongly non-reflecting at $\gs$. {\em Then\/}
  $\FP(A, \gl)$ is $<\gl$-strategically closed.
\end{lemma}

\section{The iteration}
\label{iteration-sect}

   The idea of the construction is now to define $A = \setof{\gl_i}{i < \gth}$
 (where the $\gl_i$ are as in the statement of Theorem \ref{mainthm})
 and to iterate $\FP(A \cap \gl, \gl)$ for all regular $\gl \le \gth$.
%
% Or maybe theta^+ ?
%
 A crucial point will be that the forcing at stage $\gl$ is $<\gl$-strategically
 closed, using Lemma \ref{strat} and the fact that in the iteration we have
 already arranged the required instances of non-reflection below $\gl$.

 We will do a ``Reverse Easton'' iteration, that is to say an iteration where
 direct limits are taken at strongly inaccessible limit stages and inverse limits
  are taken at other limit stages. We will refer to \cite{JEB} for details
 about this sort of iteration, and we will also follow the notation of that paper
 (in particular $\FP_\ga$ is the forcing up to stage $\ga$ and $\dot \FQ_\ga \in V^{\FP_\ga}$
 is the forcing at $\ga$).

 Formally, we will define $\dot \FQ_\ga$ to be $\{0\}$ if $\ga$ is not
 a regular cardinal, and to be $\FP(A \cap \ga, \ga)_{V^{\FP_\ga}}$
 if $\ga$ is regular. 
  We will collect some information about the iteration in the
 following lemma.

\begin{lemma} \label{portmanteau}
   Let $\FP_\ga$ and $\dot \FQ_\ga$ be as above,
 and let $\dot \FR_{\gb, \ga}$ be
 the canonical iteration in $V^{\FP_\gb}$ such that
 $\FP_\gb \ast \dot  \FR_{\gb, \ga}$ has a dense subset isomorphic to $\FP_\ga$.
 Then for all regular $\ga \le \gth$
\begin{enumerate}
\item  $\card{\FP_\ga} \le \ga$.

\item  $V^{\FP_\ga} \models GCH$, so in particular
 $V^{\FP_\ga} \models \card{\dot \FQ_\ga} = \ga$.

\item $\FP_{\ga+1}$ has the $\ga^+$-c.c. In addition, if $\ga$ is Mahlo,
 then $\FP_\ga$ has the $\ga$-c.c.

\item \label{cl4} 
 $V^{\FP_\ga} \models \dot \FQ_\ga \mbox{\ is $<\ga$-strategically closed}$.

\item \label{cl5} For  all regular $\gb < \ga$, $\FR_{\gb,\ga}$
 is $<\gb$-strategically closed in $V^{\FP_\gb}$.

\item $\FP_\ga$ preserves all cardinals and cofinalities.

 \end{enumerate}
\end{lemma}

\begin{proof}
   The proof will be by induction on $\ga$. Most of the proof is straightforward,
 using the results of Section 2 in \cite{JEB} to power the induction. The distinctive
 point here is in showing that clauses \ref{cl4} and \ref{cl5} hold at $\ga$, given that we have
 proved the lemma for regular cardinals less than $\ga$.

  By construction,
 $\FP_\ga$ forces that for every $\gl \in A$ and every regular cardinal $\gb \in [\gl, \ga)$ we have
 $\SNR(\gb, \gl)$. As we remarked after Lemma \ref{lemma1}, this implies that for every
 ordinal $\gg \in [\gl, \ga)$ we have $\SNR(\gg, \gl)$. By Lemma \ref{stratclos} this
 means that $\FP(A \cap \ga, \ga)$ is $<\ga$-strategically closed in $V^{\FP_\ga}$.

 Finally, to see that clause \ref{cl5} holds one should check that Theorem 2.5
 from \cite{JEB} is still true if ``$\gk$-closed'' is replaced by ``$<\gk$-strategically
 closed''. This is routine, the point is that a term for a strategy can be applied
 to a term for a condition to get a term for a stronger condition.

\end{proof}

We make some remarks about this construction.
\begin{enumerate}
\item Since cardinals and cofinalities are preserved, a witness to strong
 non-reflection added at some stage by some $\FP(\gs, \gl)$ will remain
 a witness at all subsequent stages.
\item At stage $\gth$ we forced with $\FP(A, \gth)$, so added witnesses
 to all the strong non-reflection that is claimed in Theorem \ref{mainthm}.
\end{enumerate}

\section{Preserving measurability}

 As we mentioned in the first section, we want to show that for each $i$
 the measurability of $\gk_i$ is preserved by the iteration $\FP_\gth$.
 It is enough to argue that $\gk_i$ is measurable in the extension
 by $\FP_{\gk_i^+ + 1}$, because the rest of the forcing is 
 $\gk_i^+$-strategically closed, so that the power set of $\gk_i$ does
 not change and a measure remains a measure.
 For brevity, we will denote
 $\gk_i$ by $\gk$ throughout this section.

  Let $G$ be $\FP_{\gk}$-generic over $V$, let $g$ be $\FP(A \cap \gk, \gk)$-generic
 over $V[G]$, and let $h$ be $\FP(A \cap \gk^+, \gk^+)$-generic over $V[G][g]$.
 Let $j: V \lra M$ be the ultrapower map arising from a normal measure $U$
 on $\gk$. We list some facts about $j$ and $M$, all of whose proofs can
 be found in \cite{SRK}.
\begin{enumerate}
\item $\crit(j) = \gk$.
\item ${}^\gk M \subseteq M$.
\item $H_{\gk^+} \subseteq M$.
\item $\gk^+ = \gk^+_M$.
\item $\gk^+ < j(\gk) < j(\gk^+) < \gk^{++}$.
\item $M = \setof{ j(F)(\gk) }{\hbox{$F \in V$ and  $\dom(F) = \gk$}}$.
\end{enumerate}

   The strategy of the proof will be to define, in $V[G][g][h]$, an extension
 of $j: V \lra M$ to a new embedding $j:V[G][g][h] \lra N \subseteq V[G][g][h]$.
 The existence of such an extension will imply that $\gk$ is still measurable
 in $V[G][g][h]$.

 We start by comparing the iterations $\FP_{\gk^+ +1}$ and $j(\FP_{\gk^+ + 1})$.
 The forcing  $j(\FP_{\gk^+ + 1})$ is an iteration defined in $M$, forcing strong
 non-reflection at cofinalities in the set $j(A) \cap j(\gk^+)$ for all $M$-regular cardinals
 up to $j(\gk^+)$. Since $A \cap \gk$ is bounded in $\gk$ and $\gl_{i+1} > \gk^+$,
 we see that 
\[
     A \cap (\gk^+ +1) = j(A) \cap (j(\gk^+) + 1) = A \cap \gk.
\]

 By the resemblance between $V$ and $M$, if we compute the iteration $j(\FP_{\gk^+ + 1})$
 up to stage $\gk^+$ we get $\FP_{\gk^+ + 1}$. We can therefore compute a generic
 extension $M[G][g][h]$ of $M$ by using the $V$-generic filters, and observing that
 $V$-generic filters are $M$-generic.

We claim that $V[G][g][h] \models {}^\gk (M[G][g][h]) \subseteq M[G][g][h]$.
Since $\FP_{\gk + 1}$ is $\gk^+$-c.c.~every canonical $\FP_{\gk+1}$-name
 for a $\gk$-sequence of ordinals is in $M$, so that easily
$V[G][g] \models {}^\gk (M[G][g]) \subseteq M[G][g]$. The forcing $\FQ_{\gk^+}$ is
 $<\gk^+$-strategically closed in $V[G][g]$, so it adds no $\gk$-sequence of
 ordinals, and we are done.

  In $M[G][g][h]$ let $\FR = \FR_{\gk^+ + 1, j(\gk)}$  be the canonical
 factor forcing to prolong $G * g * h$ to a $j(\FP_\gk)$-generic. 
 We claim that $\FR$ is $\gk^+$-strategically closed in $V[G][g][h]$.
 This follows from the fact that $\FR$ is $<\gk^{++}_M$-strategically
 closed in $M[G][g][h]$ and the fact that
 $V[G][g][h] \models {}^\gk (M[G][g][h]) \subseteq M[G][g][h]$. The point
 is that if II plays for $\gk^+$ steps in $V[G][g][h]$ using the strategy
 from $M[G][g][h]$, then every initial segment of the play is in $M[G][g][h]$,
 so that player II does not get stuck at any stage below $\gk^+$.

 The previous claim explains why we are working in $V[G][g][h]$ rather than 
 $V[G][g]$. If we truncate $j(\FP_\gk)$ at $\gk+1$ then the rest
 of the forcing will be $<\gk^+$-strategically closed in $V[G][g]$, but
 the following stage of the proof will demand $\gk^+$-strategic closure.

 Recall that $\gk^+ < j(\gk) < \gk^{++}$. In $M[G][g][h]$ the forcing $\FR$
 is $j(\gk)$-c.c.~and has size $j(\gk)$, so there are at most $j(\gk)$ maximal
 antichains in that model. In $V[G][g][h]$ let us enumerate these antichains
 as $\seq{A_\ga: \ga < \gk^+}$. Now consider a run of the game $G(\FR, \gk^+)$
 in which player I plays the following strategy; in response to $p_{2\gg}$
 player I will choose some element $q_\gg$ of $A_\gg$ such that $p_{2\gg}$
 is compatible with $q_\gg$, and then will play $p_{2\gg+1}$ which is
 some common refinement. Player II will play according to some winning strategy; after
 $\gk^+$ steps we have built a decreasing sequence of conditions which
 clearly generates an $M[G][g][h]$-generic filter $H$.

  Now we will start to extend $j$. Define $G^+ =_{def} G * g * h *H$, which will
  be $j(\FP_\gk)$-generic over $M$.
 We attempt to define $j:V[G] \lra M[G^+]$ by $j(\dot\gt^G) =_{def} j(\dot\gt)^{G^+}$.
 We  check that this is a well-defined elementary embedding,
 using the following well-known fact.

\begin{lemma}
\label{fred}
 Let $k: M \lra N$ be an elementary embedding between two transitive
 models of ZFC. Let $\FP \in M$ be some forcing, let $k(\FP) = \FQ$,
 and suppose that we have $G$ which
 is $\FP$-generic over $M$ and $H$ which is $\FQ$-generic over $N$.
 Suppose also that $k`` G \subseteq H$. Then defining $k(\dot\gt^G)
 = k(\dot\gt)^H$ for every $\dot\gt \in M^\FP$ gives a well-defined
 elementary embedding $k:M[G] \lra N[H]$, which extends $k: M \lra N$ and
 has $k(G) = H$.
\end{lemma}
 
\begin{proof} Easy, using the Truth Lemma and the elementarity
 of $k$.
\end{proof}

  By Lemma \ref{fred}, it is enough to check that $j`` G \subseteq G^+$. $G$ is generic
 for $\FP_\gk$ which was constructed as a direct limit, so every condition 
 $p$ in $G$ has support bounded in $\gk$. Since $\crit(j) = \gk$, the condition
 $j(p)$ contains the same information as $p$,  and since $G = G^+ \restriction \gk$
 we conclude that
 $j(p) \in G^+$.

  Since $G^+ \in V[G][g][h]$, we see that $M[G^+] \subseteq V[G][g][h]$.
  Also, we know that $H$ is generic for a forcing which adds no
  $\gk$-sequences of ordinals over $M[G][g][h]$, so that
  $V[G][g][h] \models {}^\gk (M[G^+]) \subseteq M[G^+]$.

  Now we aim to lift $j$ further to get a map with domain $V[G][g]$.
 In $V[G]$ the forcing $\FQ_\gk$ has cardinality $\gk$, and is
 $<\gk$-strategically closed with at most $2^\gk$ (that is $\gk^+$) many
 maximal antichains. Since $j:V[G] \lra M[G^+]$ is elementary, in $M[G^+]$ the
 forcing $\FQ_{j(\gk)}$ is $<j(\gk)$-strategically closed with at most
 $j(\gk^+)$ maximal antichains.

  Arguing as before, $\FQ_{j(\gk)}$ is $\gk^+$-strategically closed in
 $V[G][g][h]$. Since $j(\gk^+) < \gk^{++}$ we can repeat the argument
 from the construction of $H$ to build $g^+$ which is $\FQ_{j(\gk)}$-generic
 over $M[G^+]$. But it is not clear at this point that we can lift
 $j$ onto $V[G][g]$, because it may not be the case that $j`` g \subseteq g^+$.

  We will use Silver's ``master condition'' idea. Observe
 that $g \in M[G^+]$, and that $g$ is equivalent to a function
 $p$ where $\dom(p) = (A \cap \gk) \times \gk$
 and $p(\gs, -): \ga < \gk \longmapsto p(\gs, \ga)$ witnesses $SNR(\gk, \gs)$
 for each $\gs \in A \cap \gk$.

  Recall that $\FQ_{j(\gk)}$ is defined to be $\FP(j(A) \cap j(\gk), j(\gk))$.
 We claim that $p \in \FQ_{j(\gk)}$. The support condition is
 satisfied because $\gk < j(\gk)$ and (as we saw before)
 $j(A) \cap j(\gk) = j(A \cap \gk) = A \cap \gk$.
 It is enough to show that for
 each $\gs$ we have $p(\gs, -) \in \FP(\gs, j(\gk))$, which is to say that for
 all $\gd \le \gk$ of cofinality $\gs$ there is a club in $\gd$
 on which $p(\gs, -)$ is increasing. This is easy because $V$, $M$ and $M[G^+]$
 agree about cardinals and cofinalities up to $\gk^+$.

 Since $p$ is a condition in $\FQ_{j(\gk)}$, 
 when we construct $g^+$ we can arrange that $g^+ \ni p$.
 We claim that this suffices to guarantee that $j``g \subseteq g^+$.
 This follows from the observation that $p \le j(q)$ for
 every $q \in g$, which is true because $q$ has size less than
 $\gk$ and so $j(q)$ is just a copy of $q$.

 We can now build $j:V[G][g] \lra M[G^+][g^+]$, using Lemma \ref{fred}.
Before we can finish the construction, we need one piece of information
 about this embedding.
 We claim that
\[
   M[G^+][g^+] = \setof{ j(F)(\gk) }{\hbox{$F \in V[G][g]$ and  $\dom(F) = \gk$}}.
\]
 To see this let $\dot\gt^{G^+ * g^+}$ be some element of $M[G^+][g^+]$,
 where $\dot\gt$ is a $P_{j(\gk)+1}$-name in $M$. We know that $\dot\gt = j(f)(\gk)$
 for some $f \in V$, and we may as well assume that $f(\ga)$ is a $\FP_{\gk+1}$-name
 for every $\ga < \gk$. Let us define $F \in V[G]$ by $F(\ga) =_{def} f(\ga)^{G * g}$.
 Since $j(G * g) = G^+ * g^+$, we see that $j(F)(\gk) = (j(f)(\gk))^{G^+ * g^+} 
 = \dot\gt^{G^+ * g^+}$ as required.

 We will now define a filter $h^+$ on $\FQ_{j(\gk^+)}$, by setting
\[
    h^+ =_{def} \setof{q}{(\exists p \in h) \; j(p) \le q}.
\]
 It is easy to see that $h^+$ is in fact a filter, and certainly $j`` h \subseteq h^+$
 and $h^+ \in V[G][g][h]$. We claim that $h^+$ is generic. To see this, let
 $D \in M[G^+][g^+]$ be a dense subset of $\FQ_{j(\gk^+)}$. We know  $D= j(F)(\gk)$
 for some $F \in V[G][g]$. Define $E \subseteq \FQ_{\gk^+}$
by
\[
   E = \bigcap \setof{F(\ga)}{\hbox{$F(\ga)$ is a dense subset of $\FQ_{\gk^+}$}}.
\]
  $\FQ_{\gk^+}$ is $<\gk^+$-distributive, so that $E$ is dense, and clearly
 $E \in V[G][g]$. Therefore there is some $p \in E \cap h$. Certainly $j(p) \in h^+$,
 and by elementarity $D = j(F)(\gk) \supseteq j(E)$ so that $j(p) \in h^+ \cap D$.

 In conclusion, we can define $j:V[G][g][h] \lra M[G^+][g^+][h^+]$ in the model 
$V[G][g][h]$, so that $\gk$ is still measurable in  $V[G][g][h]$.

\section{The collapse}

   To save on notation, we will now denote the model $V^{\FP_\gth}$ constructed in 
 Section \ref{iteration-sect}  by $V$. In this model we have the following
 situation. For all $i$
\begin{enumerate}
\item GCH holds.
\item $\gl_i = \cf(\gl_i) < \gl_i^+ < \gth_i = \cf(\gth_i) < \gk_i$, and $\gk_i$
 is measurable.
\item $\gl_i > (\sup\setof{\gk_j}{j < i})^{++}$.
\item $\SNR(\gth, \gl_i)$ holds for every $i$.
\end{enumerate}

  We still have to get the reflection property $\Ref(\gk_i, \gth_i, <\gth_i)$
 for every $i$. We will do this by collapsing the measurable cardinals $\gk_i$,
 using an idea from Section 7 of \cite{JEB2}. We will also check that this
 collapse does not destroy the strong non-reflection.

 Let $\FS_i =_{def} Coll(\gth_i, <\gk_i)$. We define $\FS$ to  be the
 Easton product of the $\FS_i$, to be precise
$p \in \FS$ iff $p$ is a function with 
\begin{enumerate}
\item $\dom(p) \subseteq \gth$.
\item $p(i) \in \FS_i$ for all $i \in \dom(p)$.
\item If $\gs \le \gth$ is an inaccessible cardinal
 and $i < \gs \implies \gk_i < \gs$, then
 $\dom(p) \cap \gs$ is bounded in $\gs$.
\end{enumerate}
 The ordering is the natural one.

  For each $i$, the forcing $\FS$ factorises as $\FS^l_i \times \FS_i \times \FS^u_i$,
 where $\FS^l_i$ talks about the coordinates below $i$ and $\FS^u_i$ talks
 about those above. Using Easton's Lemma and the GCH, it is easy to see that
 $\FS$ collapses cardinals in the interval $[\gth_i, \gk_i)$ to $\gth_i$
 and preserves all other cardinals. In particular $\SNR(\gth, \gl_i)$
 still holds in $V^\FS$, because $\gl_i$ is still regular and there are no new points
 of cofinality $\gl_i$ (this is easy, because by our assumptions on $\gl_i$
 we have $\card{\FS^l_i} < \gl_i$, and $\FS_i \times \FS^u_i$ is $\gl_i^+$-closed).

 For the reflection, it will suffice to check that $\Ref(\gk_i, \gth_i, <\gth_i)$
 holds in $V^{\FS^l_i \times \FS_i}$, because this model agrees with $V^\FS$
 up to $\gth_{i+1}$. We will look at $V^{\FS^l_i \times \FS_i}$ in a slightly
 different way, by writing it as $(V^{\FS_i})^{\FS^l_i}$.

 Since GCH holds in $V$ and $\gk_i$ is measurable there, the results
 of \cite{JEB2} show that $\Ref(\gk_i, \gth_i, <\gth_i)$ holds
 in $V^{\FS_i}$. Of course, $\gk_i$ is now $\gth_i^+$.
 We claim that $\Ref(\gk_i, \gth_i, <\gth_i)$ still holds in 
$(V^{\FS_i})^{\FS^l_i}$.  Observe that $\card{\FS^l_i} < \gth_i < \gk_i$,
 so that if $S$ is a stationary subset of $S^{\gk_i}_{<\gth_i}$ in $(V^{\FS_i})^{\FS^l_i}$
  then there is $T \subseteq S$ stationary with $T \in V^{\FS_i}$. By the reflection
 which holds in $V^{\FS_i}$, there is $\gg \in S^{\gk_i}_{\gth_i}$ such that
 $T \cap \gg$ is stationary in $V^{\FS_i}$. Since $\card{\FS^l_i} < \gth_i$, we see that
  $T \cap \gg$ is still stationary in $(V^{\FS_i})^{\FS^l_i}$ (and of course $\gg$ still
 has cofinality $\gth_i$). 

  We have shown that $\SNR(\gth, \gl_i)$ and $\Ref(\gk_i, \gth_i, <\gth_i)$
 hold in $V^\FS$ for all $i< \gth$. This concludes the proof of
 Theorem \ref{mainthm}.

\end{document}